\documentclass[oneside,12pt]{amsart}
\usepackage[utf8]{inputenc}
\usepackage[english]{babel}
\usepackage{graphicx}
\usepackage[colorlinks=true,citecolor=cyan,linkcolor=blue]{hyperref}
\usepackage[scale=0.65]{geometry}

\newcommand{\R}{\mathbb R}
\newcommand{\N}{\mathbb N}

\newcommand{\E}{\mathbb E}
\newcommand{\cH}{\mathcal H}
\newcommand{\abs}[1]{\lvert#1\rvert}

\newtheorem{theorem}{Theorem}
\newtheorem*{theoremp}{Theorem}
\newtheorem{coro}{Corollary}
\theoremstyle{definition}
\newtheorem*{defi}{Definition}

\numberwithin{equation}{section}

\begin{document}

\title{A sunflower anti-Ramsey theorem and its applications}

\author{Leonardo Martínez-Sandoval}
\author{Miguel Raggi}
\author{Edgardo Roldán-Pensado}
\address[L. Martínez-Sandoval, E. Roldán-Pensado]{Instituto de Matemáticas, UNAM campus Juriquilla}
\email{\href{mailto:leomtz@im.unam.mx}{leomtz@im.unam.mx}, \href{mailto:e.roldan@im.unam.mx}{e.roldan@im.unam.mx}}
\address[M. Raggi]{Escuela Nacional de Estudios Superiores,
UNAM}
\email{\href{mailto:mraggi@gmail.com}{mraggi@gmail.com}}
\subjclass[2010]{}

\begin{abstract}
A $h$-sunflower in a hypergraph is a family of edges with $h$ vertices in common. We show that if we colour the edges of a complete hypergraph in such a way that any monochromatic $h$-sunflower has at most $\lambda$ petals, then it contains a large rainbow complete subhypergraph. This extends a theorem by Lefmann, Rödl and Wysocka, but this version can be applied to problems in geometry and algebra.
We also give an infinite version of the theorem.
\end{abstract}

\maketitle

\section{Introduction}

Let $K_k(V)=(V,E)$ be the complete $k$-hypergraph on $V$. In other words, let $V$ be a set and $E=\binom{V}{k}$, where $\binom{V}{k}$ represents the family of subsets of size $k$ of $V$. 

A \emph{colouring} of $K_k(V)$ is simply a function $\Delta:E \to \N$. We think of each $i\in\N$ as a colour, so $e$ has colour $i$ if $\Delta(e)=i$. If $\Delta:E \to[c]$, then we call it a \emph{$c$-colouring}. A subset $V'$ of $V$ is \emph{monochromatic} if $\Delta(\binom{V'}{k})=\{i\}$ for some $i\in[c]$, and it is \emph{rainbow} if $\Delta(e_1)\neq\Delta(e_2)$ for every pair of distinct elements $e_1,e_2\in\binom{V'}{k}$.

Ramsey's theorem \cite{Ram1930} guarantees the existence of monochromatic subsets in sufficiently large structures. There is a quote by Motzkin \cite{Pro2005} that describes this behaviour:
\begin{quote}
\it ``Complete disorder is impossible.''\\
---T. Motzkin.
\end{quote}
The precise statement of the theorem is as follows.
\begin{theoremp}[Ramsey's theorem]
For $k$, $c$ and $n_1,\dots,n_c \in \N$, there is a number $R$
such that the following holds: if $\abs{V}\ge R$ then for every $c$-colouring $\Delta$ of $K_k(V)$ there is a colour $i$ and a monochromatic subset $V'\subset V$ with $\abs{V'}=n_i$ and $\Delta(V')=\{i\}$.
\end{theoremp}

Over the years Ramsey's theorem has proven to be one of the most useful theorems in combinatorics. It has countless applications in geometry, number theory, algebra, computer science and many other fields. For more details about Ramsey Theory we refer the reader to \cite{NR1990,Gra1990,LR2014}.

Anti-Ramsey theorems appeared in 1973 \cite{ESS1973}. In contrast to the classical Ramsey theory which seeks monochromatic subsets, anti-Ramsey theory seeks \emph{rainbow} subsets.
The following quote from \cite{JNR2005} describes this behaviour:
\begin{quote}
\it ``Complete disorder is unavoidable as well.''\\
---V. Jungić, J. Ne\v{s}et\v{r}il, R. Radoi\v{c}ić.
\end{quote}
Of course, rainbow subsets of a given size may not exist for arbitrarily large hypergraphs so some extra conditions are needed. There are several types of anti-Ramsey theorems. For a survey, we refer the reader to \cite{FMO2010}.

In this note we study an anti-Ramsey theorem by Lefmann, Rödl and Wysocka \cite{LRW1996} and generalize it to Theorems \ref{thm:main} and \ref{thm:infinite} below. Theorem \ref{thm:main} has interesting applications. We use it to obtain polynomial bounds for some problems in geometry and number theory, improving results from \cite{MR2014} and \cite{GZ2013}.

In order to state our results, we need the following definition.

\begin{defi}
A \emph{$h$-sunflower} of $K_k(V)=(V,E)$ is a family $S\subset E$ consisting of edges (called \emph{petals}) whose intersection has (at least) $h$ elements.
\end{defi}

\begin{figure}
\centering
\includegraphics[scale=0.75]{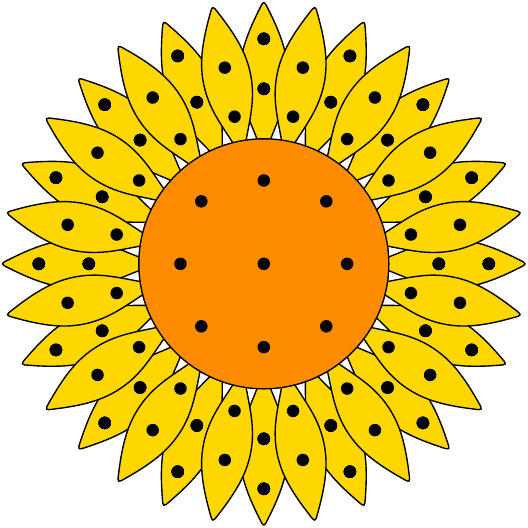}
\caption{A $9$-sunflower with $36$ petals in an $11$-hypergraph.}
\end{figure}

We are interested in colourings of $K_k(V)$ in which there are no large monochromatic $h$-sunflowers.

\begin{theorem}\label{thm:main}
Let $V$ be a set and $h < k$ and $\lambda$ be natural numbers. Let $\Delta$ be a colouring of $K_k(V)$ such that any monochromatic $h$-sunflower has at most $\lambda$ petals. Then there exists a constant $c_k$ depending only on $k$ such that if $n$ satisfies $$\abs{V}\geq c_k \lambda n^{\frac{2k-1}{k-h}}\log(n)^{-\frac{1}{k-h}},$$ there is a rainbow $V'\subset V$ with $\abs{V'}\ge n$.
\end{theorem}

The history of this theorem can be traced back to Babai \cite{Bab1985} who proved the case $k=2,h=0,\lambda=1$. Later Alon, Lefman and R\"odl \cite{ALR1992} extended Babai's result for $k\ge 2$, $h=1$ and $\lambda = 1$. Lastly Lefmann, R{\"o}dl and Wysocka \cite{LRW1996} proved it for $k\ge 2$, $h\ge 1$ and $\lambda=1$.

The proof of Theorem \ref{thm:main} is very similar to the proof of Theorem 3 in \cite{LRW1996} so we only include a sketch of it in Section \ref{sec:proofs}.

We also state a version of Theorem \ref{thm:main} for sets with infinite cardinality. It should be noted that it is not a direct consequence of its finite version, which only implies the existence of arbitrarily large rainbow subsets. The proof of this theorem can also be found in Section \ref{sec:proofs}.

\begin{theorem}\label{thm:infinite}
Let $V$ be an infinite set and let $\Delta$ be a colouring of $K_k(V)$ with no monochromatic $h$-sunflower with $\lambda$ petals. Then there is an infinite rainbow subset of $V$.
\end{theorem}

Note that the condition on $\lambda$ cannot be relaxed to only forbid finite monochromatic sunflowers. It is easy to construct an example showing this: for any edge $e$ of $K_k(\N)$, define $\Delta(e)=\max(e)$.

There are some interesting applications of Theorem \ref{thm:main}. These are stated below and their proofs are included in Section \ref{sec:applications}.

We begin with applications to geometry, specifically to Erd\H os-Szekerez type problems \cite{ES1935}. Here the objective is to find, within a large set of points in general position, a subset of a given size with specific properties. Usually the existence of such a subset can be shown using Ramsey's theorem but the bounds on the size of the original set are usually extremely large. The first corollary gives polynomial bounds to three of these problems.

\begin{coro}\label{coro:geometry}
Let $d$ and $n$ be a positive integers. There exists a constant $c_d$ depending only on $d$ such that if $X$ is a set of points in $\R^d$ for which
\[
\abs{X}\geq c_d \frac{n^{2d+1}}{\log(n)}
\]
then the following hold:
\begin{enumerate}
\item If no $d+2$ points of $X$ lie on a $(d-1)$-sphere, then there exists $Y\subset X$ with $n$ points such that the circumradii of all the $d$-simplices defined by elements of $Y$ are distinct.
\item If no $d+1$ points of $X$ lie on a hyperplane then there exists $Y\subset X$ with $n$ points such that the volumes of all the $d$-simplices defined by elements of $Y$ are distinct.
\item If no $d+1$ points of $X$ lie on a hyperplane then there exists $Y\subset X$ with $n$ points such that the similarity types of all the $d$-simplices defined by elements of $Y$ are distinct. 
\end{enumerate}
\end{coro}

Part (1) of this theorem improves the main result in \cite{MR2014}, in which the bound $O(n^9)$ was given for the case $d=2$. Part (2), with $d=2$, improves a bound given in \cite{GZ2013} that uses Ramsey's theorem taking it down from $O(e^{e^{n^{186}}})$ to $O(n^5/\log(n))$. Furthermore, this version generalises to any number of dimensions.

As a second application we show a result regarding polynomials.

\begin{coro}\label{coro:polynomial}
Let $\mathbb F$ be a field and $p\in\mathbb F[x,y]$ be a symmetric polynomial of degree $d\ge 1$. There exists a constant $c_d$ depending only on $d$ such that if $X\subset\mathbb F$ and
\[\abs{X}\geq c_d \frac{n^3}{\log(n)},\]
then we can choose a subset $Y\subset X$ of size $n$ such that all the evaluations of $p$ in $2$ distinct elements of $Y$ are distinct.
\end{coro}

If the polynomial $p$ has degree $1$ then $p$ is simply the sum of $2$ variables. An analogous theorem for sums of $k$ distinct elements in Abelian groups also holds and has a similar proof, although this was already a consequence of the theorem in \cite{LRW1996}.

Finally we consider $B_2$-sequences, also called Sidon sequences (see \cite{Sid1932,ET1941}). These are increasing sequences (possibly finite) of positive integers $a_1<a_2<\dots$ such that the sums $a_i+a_j$ are distinct for every $i\le j$. An equivalent definition is that the differences $a_j-a_i$ are distinct for every $i\le j$.

\begin{coro}\label{coro:sidon}
There is a constant $c>0$ such that any set $X$ of positive integers with $\abs{X}\ge c n^3/\log(n)$ contains a $B_2$-sequence of length $n$.
\end{coro}

For comparison, it is known that for every $\varepsilon>0$ there is a large enough $n$ such that the set $\{1,2,\dots,(\frac{1}{2}+\varepsilon)n^2\}$ contains a $B_2$-sequence of length $n$.

\section{Proofs of Theorems \ref{thm:main} and \ref{thm:infinite}}\label{sec:proofs}

As stated before, the proof of Theorem \ref{thm:main} is very similar to the proof Theorem 3 in \cite{LRW1996} and only needs minor adjustments. Therefore, we only include a sketch of the main idea. The proof uses a theorem by Ajtai et al. \cite{AKP+1982} that gives a lower bound on the independence number of an \emph{uncrowded} hypergraph. An uncrowded hypergraph is one which has no cycles of length $2$, $3$ or $4$.

\begin{proof}[\textbf{Sketch of the proof of Theorem \ref{thm:main}}]
Let $\Delta:K_k(V)\to\N$ be a colouring with no monochromatic $h$-sunflowers with more than $\lambda$ petals. We wish to find a large rainbow subset of $V$. For convenience, denote $N=\abs{V}$.

We construct a $2k$-hypergraph $\cH$ as follows: The vertices are $V$ and $e \subset \binom{V}{2k}$ is an edge if and only if it contains two different sets $A,B \in \binom{V}{k}$ with $\Delta(A)=\Delta(B)$. Then an independent set on $\cH$ corresponds to a rainbow set of $K_k(V)$.

There is a simple probabilistic argument to get rid of all the edges of $\cH$ while keeping a good portion of the vertices, but this method yields the desired bound without the $\log(n)^{\frac{1}{k-h}}$ factor. We chose to illustrate it here, as the proof of the tighter bound involves applying the same method to instead get rid of all 2,3,4-cycles of $\cH$ first while keeping a larger amount of vertices, and then using the lower bound on the independence number of an uncrowded hypergraph by Ajtai et al.

The sunflower condition implies that the number of edges of $\cH$ can be bounded as follows:
$$|E(\cH)| \leq \binom{N}{k}\binom{k}{h}\lambda\binom{N-2k+h}{h} \approx c_k\lambda N^{k+h},$$ 
but these edges have size $2k$, which is larger than $k+h$, since $h \leq k-1$.

We use the probabilistic method: choose $p\in(0,1)$ and then choose each vertex independently at random with probability $p$ (erase with probability $1-p$). Since $2k > k+h$, if we choose $p$ adequately we can ``destroy'' all but a few edges. Let $V', E'$ denote the random variables which count the number of vertices and edges respectively.

By linearity of the expected value $\E$, we have:
$$\E(V') = pN,\qquad \E(E') = p^{2k}E \leq p^{2k}N^{k+h}.$$

After this process we would like to obtain a hypergraph in which the number of edges is in $o(pN)$, so we can delete the remaining edges ``by hand'' (by, say, deleting one vertex from each). This can be assured if the expected value of $E'$ is in $o(pN)$. Equality between $p^{2k}N^{k+h}$ and $pN$ is achieved when $p = N^{-\frac{k+h-1}{2k-1}}$, so that $pN = N^{\frac{k-h}{2k-1}}$. Choosing a slightly smaller $p$ achieves the desired result. With some technical details omitted, we can get an independent set of $\cH$ of size in $\Theta(pN)$, which is a rainbow of the original. 

The gritty details for the improved bound can be found in \cite{LRW1996}, and their proof works almost word for word for this generalisation, the only difference being that some expressions have to be multiplied by the constant $\lambda$, which does not affect the asymptotic behaviour.
\end{proof}

Now we prove the generalisation of Theorem \ref{thm:main} to infinite hypergraphs. This theorem is implied by a theorem of Erd\H os and Rado \cite{ER1950}.

\begin{theorem}[Erd\H os-Rado theorem]\label{thm:ER}
Let $\Delta:\binom{\N}{k}\to\N$ be any colouring of the $k$-sets of $\N$. Then there exists an infinite subset $X\subset\N$ and a set of indices $I=\{i_1,i_2,\dots,i_r\}$ with $1\le i_1<i_2<\dots<i_r\leq k$ such that two ordered subsets $A,A'\in \binom{X}{k}$ have the same colour if and only if $a_{i_j} = a_{i_j}'$ for every $j\in[r]$.
\end{theorem}
 
For example, when $I=\emptyset$, we have a monochromatic structure because the colour of a $k$-set does not depend on any of the elements, and when $I=[k]$ we have a rainbow structure because the colour of a $k$-set depends on every element.
 
Now we are ready to prove Theorem \ref{thm:infinite}.

\begin{proof}[\textbf{Proof of Theorem \ref{thm:infinite}}]
Without loss of generality we may assume $h=k-1$. Suppose $\Delta:\binom{\N}{k} \to \N$ has the property that any monochromatic $(k-1)$-sunflower has at most $\lambda$ petals. Let $X\subset \N$ and $I\subset[k]$ be sets as in Theorem \ref{thm:ER}, we wish to prove that $I=[k]$. For sake of contradiction, assume that some number $a\in[k]$ is not in $I$. Fix numbers $n_1<n_2<\dots<n_{k-1}$ in $X$ such that there are distinct $m_1,m_2,\dots,m_{\lambda+1}\in X$ strictly between $n_{a-1}$ and $n_a$. Then the $k$-sets $\{n_1,n_2,\dots,n_{k-1}\}\cup\{m_i\}$ for $i=1,\dots,\lambda+1$ form a monochromatic $(k-1)-$sunflower with $\lambda+1$ petals.
\end{proof}

Theorem \ref{thm:infinite} makes it clear that the corollaries above have corresponding infinite versions.

\section{Proofs of corollaries}\label{sec:applications}

In this section we will only use the case $h=k-1$ of Theorem \ref{thm:main}, so the condition on the size of $V$ becomes
$$\abs{V}\geq c_k \frac{\lambda n^{2k-1}}{\log(n)}.$$
The proofs are straightforward, all that needs to be done is define a colouring of the $k$-tuples of a set and show that the monochromatic sunflower condition is satisfied for some $\lambda$. Then Theorem \ref{thm:main} gives us the desired subset.

\begin{proof}[\textbf{Proof of Corollary \ref{coro:geometry}}]
For the three parts we tale $k=d+1$ but colour the $(d+1)$-tuples in different ways.

For part (1), we colour a $(d+1)$-tuple by its circumradius. Given $d$ points in $\R^d$, there are at most $2$ spheres of a given radius which pass through them, so we take $\lambda=2$. 

In part (2), we colour a $(d+1)$-tuple by its volume. Given $d$ points, the set of points that complete this $d$-tuple into a simplex of a prescribed volume is the union of two hyperplanes, since there are at most $d$ points on each, we set $\lambda=2d$.

For the last part, we colour a $(d+1)$-tuple by its similarity type. Fix $x_1,\dots,x_d\in X$ and let $S$ and $S'$ be sets of $d+1$ points that determine a certain similarity type but with different orientations. If $x_{d+1}$ makes $\{x_1,\dots,x_{d+1}\}$ have the same similarity type as $S$, then it determines a unique injective function $\{x_1,\dots,x_d\}\to S$ or a unique injective function $\{x_1,\dots,x_d\}\to S'$. Since there are $2(d+1)!$ of these functions, we may take $\lambda=2(d+1)!$.
\end{proof}

\begin{proof}[\textbf{Proof of Corollary \ref{coro:polynomial}}]
Since $p$ is symmetric we can colour a pair of distinct elements $x,y\in Y$ by the value of $p(x,y)$. Now write $p(x,y)=\sum_{i=0}^d q_i(y)x^i$. Let $j\ge 1$ be the smallest index for which $q_j(y)\neq 0$, $Z$ be the zero-set of $q_j(y)$, and $Y=X\setminus Z$. Note that, because the degree of $p$ is $d$, $Z$ has at most $d$ elements. For every fixed $y_0\in Y$, $p(x,y_0)$ is a polynomial on $x$ with degree between $1$ and $d$, therefore the equation $p(x,y_0)=a$ has at most $d$ solutions. Thus, we can set $\lambda=d$.
\end{proof}

\begin{proof}[\textbf{Proof of Corollary \ref{coro:sidon}}]
We colour a pair $\{x,y\}\subset X$ by $\abs{x-y}$. For each $y\in X$ there are at most two elements $x\in X$ that produce the same value of $\abs{x-y}$, so we can take $\lambda=2$.
\end{proof}

\end{document}